\documentclass[reqno]{amsart}

\usepackage{amsmath,amsthm,amsfonts,amssymb,young}
\usepackage{array, boldline, makecell, booktabs}
\usepackage{bm}
\usepackage{ytableau}
\usepackage{euscript, mathrsfs} 
\usepackage{tikz}
\usetikzlibrary{backgrounds}
\usepackage{amsmath}
\usepackage{mathtools}
\usepackage[colorlinks]{hyperref}
\usepackage{cleveref}
\usepackage[margin=1in]{geometry}
\numberwithin{equation}{section}    
\usepackage[all]{xy}
\usepackage{graphicx,import}
\usepackage{amscd}
\usepackage{color}
\usepackage{enumerate}
\usepackage{fourier}
\usepackage{tikz}
\usepackage{cases}
\usepackage{cite}
\usepackage[T1]{fontenc}
\usepackage{tikz-cd}
\usepackage{leftindex}
\numberwithin{equation}{section}

\newcommand{\LL}{{\mathrm{L}}}
\newcommand{\KL}{{\mathrm{KL}}}

\newtheorem{thm}{Theorem}[section]

\theoremstyle{definition}
\newtheorem{example}[thm]{Example}
\newtheorem{definition}[thm]{Definition}
\newtheorem{conj}[thm]{Conjecture}
\newtheorem{rmk}[thm]{Remark}

\title{Combinatorial description of Lusztig $q$-weight multiplicity}
\author{Seung Jin Lee}
\date{August 2024}
\address{
Department of Mathematical Sciences, Research Institute of Mathematics, 
Seoul National University, Seoul 151-747,
South Korea}
\email{lsjin@snu.ac.kr}


\begin{document}

\begin{abstract}
We conjecture a precise relationship between Lusztig $q$-weight multiplicities for type $C$ and Kirillov-Reshetikhin crystals. We also define $\mathfrak{gl}_n$-version of $q$-weight multiplicity for type $C$ and conjecture the positivity.
\end{abstract}
\maketitle

\section{Introduction}

The Kostka-Foulkes polynomials $K_{\lambda\mu}(t)$ are coefficients for the modified Hall-Littlewood polynomials when expressed in the Schur basis. They are primarily studied in algebraic combinatorics and representation theory. In 1978, Lascoux and Schützenberger \cite{LS78} showed that 
$$K_{\lambda\mu}(t)=\sum_{T\in \textrm{SSYT}(\lambda,\mu)}t^{\textrm{charge}(T)},$$ where the sum is taken over all semistandard Young tableaux of shape $\lambda$ and weight $\mu$. Since the discovery of the charge formula \cite{LS78, LS81}, there have been various interpretations of the charge statistic.

Kostka-Foulkes polynomials can be generalized in mainly two ways. One important generalization is the Macdonald-Kostka polynomials $K_{\lambda\mu}(q, t)$, which are the coefficients for the modified Macdonald polynomials expressed in terms of Schur functions. When $q=0$, the Macdonald-Kostka polynomial reduces to the Kostka-Foulkes polynomial. Regarding positivity, Haiman \cite{Hai01} showed that $K_{\lambda\mu}(q, t)$ lies in $\mathbb{Z}_{\geq 0}[q, t]$ by investigating Hilbert schemes. Macdonald polynomials can also be expanded positively in terms of LLT polynomials \cite{HHL05}, and Haiman and Grojnowski \cite{HG06} demonstrated that LLT polynomials are Schur positive by proving that the coefficients appear in the theory of Hecke algebras and Kazhdan-Lusztig polynomials, which are known to be nonnegative. Despite various proofs of the positivity of Kostka-Foulkes polynomials, no combinatorial (manifestly positive) formula is currently known.

Another way to generalize the Kostka polynomials is to consider those for other types, known as Lusztig's $q$-weight multiplicities \cite{Lus83}, which are the main focus of this paper. For a simple Lie algebra $\mathfrak{g}$ and dominant weights $\lambda$ and $\mu$, the Lusztig's $q$-weight multiplicity $\KL^{\mathfrak{g}}_{\lambda\mu}(q)$ is defined by the formula \begin{align}\label{Def_KL} D_{w_0}\left(e^\lambda \prod_{\alpha \in \Delta^+} \frac{1}{1 - q e^\alpha}\right) = \sum_\mu \KL^{\mathfrak{g}}_{\lambda\mu}(q) \chi^\mu, \end{align} where $w_0$ is the longest element in the Weyl group of $\mathfrak{g}$, $\Delta^+$ is the set of positive roots, $D_w$ is the Demazure operator, and $\chi^\mu$ is the irreducible character for $\mathfrak{g}$ indexed by $\mu$.  Note that $D_{w_0}(e^\lambda) = \chi^\lambda$ is the Weyl character formula. Lusztig's $q$-weight multiplicities are known to be nonnegative due to the theory of affine Kazhdan-Lusztig polynomials, but manifestly positive formulas are not generally known. For type $C$, Lecouvey \cite{Lec05} constructed the cyclage on Kashiwara-Nakajima tableaux \cite{KN94} and conjectured that the statistics defined by the cyclage might provide a manifestly positive formula for Lusztig's $q$-weight multiplicity in type $C$. Additionally, Lecouvey and Lenart \cite{LL20} found a combinatorial formula for $\KL^{C_n}_{\lambda\mu}(q)$ when $\mu = 0$.

In this paper, we describe the precise relationship between Lusztig's $q$-weight multiplicities for type $C_n$ and the energy functions in the affine KR crystals for affine type $B^{(1)}_g$ for large $g$. Note that Nakayashiki and Yamada \cite{NY97} demonstrated that the charge for type $A$ corresponds to the energy function in the affine KR crystals for affine type $A$. To describe the conjecture, we introduce semistandard oscillating tableaux (SSOT) and their realization as certain classical highest weights in the tensor products of Kirillov-Reshetikhin column crystals.

Let $\overline{\mathfrak{g}}$ be an affine Lie algebra and $U'_q(\overline{\mathfrak{g}})$ the corresponding quantum algebra without derivation. Kirillov-Reshetikhin crystals (KR crystals) are crystal bases $B^{r,s}$ \cite{Kas91,Lus90} of a certain subset of irreducible finite-dimensional $U'_q(\overline{\mathfrak{g}})$-modules known as Kirillov-Reshetikhin (KR) modules $W^{r,s}$ \cite{HKOTY98,HKOTY02}. For a partition $\gamma$, we define the tensor products of KR column crystals $B^t_\gamma$ as $B^{\gamma_1,1} \otimes \cdots \otimes B^{\gamma_n,1}$. For each element $r$ in $B^t_\gamma$, we can associate a non-negative integer $E(r)$, called the energy, which is invariant under the Kashiwara operators $e_i$ and $f_i$ for $i \neq 0$. In our specific case, $\overline{\mathfrak{g}}$ is of type $B^{(1)}_g$ for large $g$, but it will be evident that the idea can be generalized to other types (see Section 3 for the stable version).

On the other hand, the author \cite{Lee23} constructed a bijection between King tableaux of shape $\lambda$ and weight $\mu$ and SSOTs of shape $\tilde{\lambda}$ and $\tilde{\mu}$ with at most $g$ columns, where $\tilde{\lambda}$ is the rectangular complement partition in $(g^n)$ for large $g$. Note that the number of such objects is equal to the weight multiplicity of weight $\mu$ in the irreducible $\mathfrak{sp}_n$-module indexed by $\lambda$.

The goal of this paper is to define a certain statistic $E(T)$, also called energy, for each SSOT $T$ by making an injective map from the set of SSOT's of shape $\tilde{\lambda}$ and $\tilde{\mu}$ to the set of some classical highest weights in $B^t_{\tilde{\mu}}$ of weight $\lambda$ and then take the energy function defined on $B^t_{\tilde{\mu}}$. To be more precise, an SSOT consists of oscillating horizontal strip $(\alpha,\beta,\gamma)$ where
$\alpha, \beta, \gamma$ are partitions such that $\alpha$ and $\gamma$ are contained in $\beta$, $\beta/\alpha$ and $\beta/\gamma$ are horizontal strips, and the size of a oscillating horizontal strip is $\ell=2|\beta|-|\alpha|-|\gamma|$ where $|\alpha|$ is defined by the sum of parts of $\alpha$. For each oscillating horizontal strips $(\alpha,\beta,\gamma)$, we associate an element in $B^{\ell,1}$ which has an entry $i$ if and only if there is a cell at $i$-th column of $\beta/\alpha$, and has an entry $\overline{i}$ if and only if there is a cell at $i$-th column of $\beta/\gamma$. For example, for an oscillating horizontal strip $((1),(3,1),(2,1))$ we associate an element 
$$\begin{ytableau}
    1\\2\\3\\ \overline{3} 
\end{ytableau}=\begin{ytableau}
    1\\2 
\end{ytableau}$$
 in $B^{4,1}$. The right-hand side uses the sage notation, making the column admissible. 

Let $\epsilon^C(T)$ be the maximum number of columns among partitions in the SSOT $T$. In KR crystal language, $\epsilon^C(T)$ is the maximum number $i$ appearing in the element $T$ in $B^t_{\tilde\mu}$. Then we have the following conjecture:

\begin{conj}\label{conj1}
$$\KL^{C_n}_{\lambda,\mu}(q)=\sum_{T} q^{E(T)}$$
where the sum runs over all SSOT $T$ of shape $\tilde{\lambda}$ and weight $\tilde{\mu}$ and $\epsilon^C(T)\leq g$.
\end{conj}

In general, defining and computing the energy function can be challenging, so we provide an example where the energy function is easy to compute, which is the case $\tilde\mu=(1^n)$.
For an element 
$$    
a_n
\otimes 
a_{n-1}
\otimes\cdots\otimes
a_1
$$
in $\left(B^{1,1}\right)^{\otimes n}$, written with the sage notation, we define the energy function by $\sum_{i=1}^{n-1} (n-i)H(a_{i+1},a_i)$ where

$$H(b,a)=\begin{cases}
    2 &\textrm{ if } a=1 \textrm{ and } b=\overline{1}\\
    1 &\textrm{ if } b \succeq a \textrm{ and } (b,a)\neq (\overline{1},1) \\
    0 &\textrm{ if } b \prec a \\
    
\end{cases}$$
under the order
$$1 \prec 2 \prec \cdots n \prec \overline{n} \prec \cdots \prec\overline{1}.$$
Note that when all $a_i$'s are positive, the energy function is the same as the charge statistic of a standard young tableau in Lascoux-Schützenberger formula. In general, Conjecture \ref{conj1} generalizes Lascoux-Schützenberger's charge formula for semistandard Young tableaux, which applies to the case where $|\lambda| = |\mu|$. This conjecture is a significant generalization of the charge formula for type $A$, as it not only describes the $q$-weight multiplicities for type $C$ but also establishes a connection with affine crystal theory.
Also note that Conjecture \ref{conj1} naturally implies the following monotonicity:
$$\KL_{\lambda+(1^n),\mu+(1^n)}^{C_n}(q)-\KL_{\lambda,\mu}^{C_n}(q) \in \mathbb{Z}_{\geq 0}[q].$$

The following is an easy example of Conjecture \ref{conj1}
\begin{example}

    Let $n=3, \lambda=(1,1)$ and $\mu=\emptyset$. When $g=1$, $\tilde\lambda=(1)$ and $\tilde\mu=(1,1,1)$ and all classcal highest weights in $\left( B^{1,1} \right)^{\otimes 3}$ with the weight $\tilde\lambda$ are

    \begin{align*}
        \overline{1}\otimes 1 \otimes 1, \quad 1 \otimes \overline{1} \otimes 1, \quad \overline{2}\otimes 2\otimes 1.
    \end{align*}
    Since the values of $\epsilon^C$ for the first two elements are $1$ and the value of $\epsilon^C$ for the third element is $2$, with corresponding energies being 2, 4, and 3, respectively, we have the following:
    
    $$\KL_{(1,1),\emptyset}^{C_n}(q)=q^2+q^4.$$
    Similarly, we have
    $$\KL_{(2,2,1),(1,1,1)}^{C_n}(q)=q^2+q^3+q^4.$$
\end{example}
In Section 3, we prove the conjecture for the stable case by relating $q$-weight multiplicities to a one-dimensional sum and the $X  = \leftindex^{\infty}{\KL}$ theorem by Shimozono and Lecouvey \cite{LS07}. The proof in Section 3 applies not only to type $C$ but also to all nonexceptional types.

During the course of this work, the author identified another natural version of the $q$-weight multiplicity, which we term the $\mathfrak{gl}_n$ $q$-weight multiplicities for nonexceptional types. Consider a function $L_A$ on the set of positive roots, where $L_A(\alpha) = 1$ if $\alpha$ is a positive root of type $A$, namely $\epsilon_i - \epsilon_j$ for $i < j$, and $0$ otherwise. Then the $\mathfrak{gl}_n$ $q$-weight multiplicity for $\lambda,\mu$ of type $\mathfrak{g}$ is defined by the following:

\begin{definition}
For a function $L$ from the set of positive roots to $\mathbb{Q}$, define
$$ D_{w_0}\left(e^\lambda \prod_{\alpha \in \Delta^+} \frac{1}{1 - q^{L(\alpha)} e^\alpha}\right) := \sum_\mu \KL^{\mathfrak{g}, L}_{\lambda\mu}(q) \chi^\mu.
$$
Then $\mathfrak{gl}_n$ $q$-weight multiplicity for $\lambda,\mu$ of type $\mathfrak{g}$ is defined by $\KL_{\lambda\mu}^{\mathfrak{g},L_A}$
\end{definition}
\begin{conj}\label{conj2}
$\KL^{C_n,\LL_A}_{\lambda,\mu}(q)$ is in $\mathbb{Z}_{\geq 0}[q]$.

\end{conj}

Proofs of both Conjecture \ref{conj1} and Conjecture \ref{conj2} will be presented in a separate paper \cite{CKL+}. In this paper, we investigate the connection between Conjectures \ref{conj1} and \ref{conj2} and known results, particularly in relation to affine crystal theory and $X=\leftindex^{\infty}{\KL}$ theorem. 

\section*{acknowledgement}
The author formulated Conjecture \ref{conj1} around December 2019 and shared with Mark Shimozono and Jaehoon Kwon. I thank them for keeping our discussions confidential and for their valuable insights. I also appreciate Hyunjae Choi, Donghyun Kim for helpful discussions. This work is supported by the National Research Foundation of Korea (NRF) grant funded by the Korean government (MSIT) (No.0450-20240021)

\section{Lusztig $q$-weight mulitiplicity and 1-dimensional sum}
Let $\overline{\mathfrak{g}}$ be an affine algebra for nonexceptional type. In this section, we explain $X=\leftindex^{\infty}{\KL}$ theorem briefly and relate the theorem with Conjecture \ref{conj1} and \ref{conj2}.\\

First, note that Equation (\ref{Def_KL}) is equivalent to the follwing formula:
$$\KL^{\mathfrak{g},L}_{\lambda\mu}(q)= \sum_{w\in W} (-1)^w [e^{w(\lambda+\rho)-\mu-\rho}]\prod_{\alpha\in\Delta^+} \frac{1}{1-q^{L(\alpha)}e^\alpha}$$
where $[e^\beta]f$ denotes the coefficient of $e^\beta$ in $f \in \mathbb{Z}[P]$. The stable version is defined by
$$\leftindex^{\infty}{\KL}^{\mathfrak{g},L}_{\lambda\mu}(q)= \sum_{w\in S_n} (-1)^w [e^{w(\lambda+\rho)-\mu-\rho}]\prod_{\alpha\in\Delta^+} \frac{1}{1-q^{L(\alpha)}e^\alpha},$$
where $S_n$ is the symmetric group, as a subgroup of $W$.
The reason why this definition is called stable version is because for large $k$, we have
$$\leftindex^{\infty}{\KL}^{\mathfrak{g},L}_{\lambda\mu}(q)=\KL^{\mathfrak{g},L}_{\lambda+(k^n),\mu+(k^n)}(q).$$

One-dimensional (1-d) sums $X$ are graded tensor product multiplicities for affine Kac-Moody algebras, which arise from two-dimensional solvable lattice models \cite{KKMMNN92}, and which may be defined using the combinatorics of affine crystal graphs \cite{HKOTY98,HKOTY02}. For any nonexceptional family of affine algebras, the 1-d sums have a large rank limit which is called the stable 1-d sums. There are only four distinct kinds of stable 1-d sums \cite{Shi05,SZ06}, and they are labeled by the four partitions $\lozenge \in \{\emptyset,(1),(2),(1,1)\}$ having at most $2$ cells. See \cite{LS07} for more description. To explain Conjecture \ref{conj1} we only need the case $\lozenge=(2)$ which corresponds to the case when $\overline{\mathfrak{g}}$ is of type $B^{(1)}_g , A^{(2)}_{2g-1}$ or $D^{(1)}_g$ for large $g$. To be precise, let 
$$B_\mu=B^{1,\mu_1}\otimes \cdots \otimes B^{1,\mu_\ell}$$
$$B^t_\mu=B^{\mu_1,1}\otimes \cdots \otimes B^{\mu_\ell,1}$$
where $\ell=\ell(\mu)$.

For any $\lozenge$ and for any tensor product $B$ of KR modules, the $X$ polynomial is defined by
$$\overline{X}^\lozenge_{\lambda,B}= \sum q^{E(T)}$$
where $T$ runs over all classical highest elements in $B$ for type $\lozenge$, of weight $\lambda$. In this paper, $B$ is either $B_\mu$ or $B^t_\mu$.

The following is $X=\leftindex^{\infty}{\KL}$ theorem in \cite{LS07}.

\begin{thm}\label{thm_X=KL} Let $||\mu||=\sum_{i=1}^n (i-1)\mu_i$. Then we have
\begin{align*}\overline{X}^\lozenge_{\lambda,B_\mu}(q)=q^{||\mu||+|\mu|-|\lambda|}\leftindex^{\infty}{\KL}^{\mathfrak{g}_n}_{\tilde\lambda \tilde\mu}(q^{-1})\end{align*}
\end{thm}

\section{proof of Conjecture \ref{conj1} and \ref{conj2} for the stable case}

There is a follwing duality in $X$.
\begin{thm}\label{thm_Xdual}\cite[Theorem 10.10]{LOS12}
$$\overline{X}^{\lozenge^t}_{\lambda^t,B^t_\mu}(q)=q^{||\mu||+|\mu|-|\lambda|} \overline{X}^\lozenge_{\lambda,B_\mu}(q^{-1})$$
\end{thm}
Therefore, we have 
$$\leftindex^{\infty}{\KL}^{\mathfrak{g}_n}_{\lambda \mu}(q)=\overline{X}^{(1,1)}_{\tilde\lambda,B^t_{\tilde\mu}}(q)$$
by Theorem \ref{thm_X=KL} and \ref{thm_Xdual}, and Conjecture \ref{conj1} for the stable case follows.

To prove Conjecture \ref{conj2} for the stable case, consider the automorphism $\phi(e^{\epsilon_i})=q^{1/2} e^{\epsilon_i}$. Note that the automorphism $\phi$ commutes with the action of $S_n$, but not with the action of $W$. For type $C$ positive roots, $\phi(e^{\alpha})$ is $e^\alpha$ if $\alpha$ is of type $A$, or $q e^\alpha$ otherwise. Therefore, we have $$ \phi(e^{w(\lambda)})=\phi(e^{\lambda})= q^\frac{\langle \lambda,\omega_n\rangle}{2} e^\lambda$$
where $w\in S_n$ and $\omega_n=(1^n)$. Also, we have
$$\phi\left(\prod_{\alpha\in\Delta^+} \frac{1}{1-q^{L_A(\alpha)}e^\alpha}\right)=\prod_{\alpha\in\Delta^+} \frac{1}{1-qe^\alpha}.$$
Therefore, we have
\begin{thm} \label{thm_12}Assume that $\mathfrak{g}$ is nonexceptional. Then 
\begin{align*}q^{|\lambda|-|\mu|}\leftindex^{\infty}{\KL}_{\lambda,\mu}^{\mathfrak{g},L_A}(q)=\leftindex^{\infty}{\KL}_{\lambda,\mu}^{\mathfrak{g}}(q).\end{align*}
\end{thm}
Therefore, Conjecture \ref{conj2} for the stable case holds.
\begin{rmk}
There is no such an elegant identity in \ref{thm_12} for non-stable cases, because $\phi(e^{w(\alpha)})$ is not equal to $\phi(e^\alpha)$ for $w\in W$ in general.
\end{rmk}
\section{Further directions}
One natural question is whether Conjectures \ref{conj1} and \ref{conj2} can be generalized to other types. Partial progress on this topic will appear in \cite{CKL+}. Applications of Conjecture \ref{conj1} include its connections with parabolic $q$-weight multiplicities \cite{LOS12}, type $C$ analogues \cite{LS+} of Catalan functions \cite{BMPS198}, and potentially with Macdonald polynomials and LLT polynomials for type $C$.

Another question is whether there is a geometric or representation-theoretic interpretation of the duality that appears in Conjecture \ref{conj1}. It would be interesting to explain Theorem \ref{thm_Xdual} without relying on the $X = K$ result, where $K$ denotes Kostka polynomials. The proof of Theorem \ref{thm_Xdual} utilizes the $X = K$ theorem and the duality for Kostka polynomials.

Lastly, Conjecture \ref{conj1} is not compatible with Lecouvey's conjecture in \cite{Lec05}, but it seems to be compatible with the work \cite{LL20} by Lecouvey and Lenart when $\mu=0$. It would be interesting to generalize the computation of the energy in their way.

\bibliographystyle{alpha}  
\bibliography{main.bib}

\newcommand{\etalchar}[1]{$^{#1}$}
\begin{thebibliography}{KKM{\etalchar{+}}92}

\bibitem[BMPS19]{BMPS198}
Jonah Blasiak, Jennifer Morse, Anna Pun, and Daniel Summers.
\newblock Catalan functions and $k$-schur positivity.
\newblock {\em J. of Amer. Math. Soc.}, 32:921--963, 2019.

\bibitem[CKL24]{CKL+}
Hyeonjae Choi, Donghyun Kim, and Seung~Jin Lee.
\newblock Lusztig $q$-weight multiplicity and affine crystals.
\newblock {\em In preparation}, 2024+.

\bibitem[Hai01]{Hai01}
Mark Haiman.
\newblock Hilbert schemes, polygraphs and the macdonald positivity conjecture.
\newblock {\em J. Amer. Math. Soc.}, 14:941–--1007, 2001.

\bibitem[HG06]{HG06}
Mark Haiman and Ian Grojnowski.
\newblock Affine hecke algebras and positivity of llt and macdonald polynomials.
\newblock {\em preprint}, 2006.

\bibitem[HHL05]{HHL05}
James Haglund, Mark~D. Haiman, and Nicholas~A. Loehr.
\newblock A combinatorial formula for macdonald polynomials.
\newblock {\em J. Amer. Math. Soc.}, 18(03):735--762, 2005.

\bibitem[HKO{\etalchar{+}}98]{HKOTY98}
G.~Hatayama, A.~Kuniba, M.~Okado, T.~Takagi, and Y.~Yamada.
\newblock Remarks on fermionic formula.
\newblock 1998.

\bibitem[HKO{\etalchar{+}}02]{HKOTY02}
G.~Hatayama, A.~Kuniba, M.~Okado, T.~Takagi, and Y.~Yamada.
\newblock Affine crystals and vertex models.
\newblock {\em Progress in Mathematical Physics}, 23, 2002.

\bibitem[Kas91]{Kas91}
Masaki Kashiwara.
\newblock On crystal bases of the $q$-analogue of universal enveloping algebras.
\newblock {\em Duke Math.}, 63:465--516, 1991.

\bibitem[KKM{\etalchar{+}}92]{KKMMNN92}
Seok-Jin Kang, Masaki Kashiwara, Kailash~C. Misra, Tetsuji Miwa, Toshiki Nakashima, and Atsushi Nakayashiki.
\newblock Affine crystals and vertex models.
\newblock {\em International J. of Modern Physics}, 7(Suppl. 1A):449–--484, 1992.

\bibitem[KN94]{KN94}
Masaki Kashiwara and Toshiki Nakashima.
\newblock Crystal graphs for representations of the $q$-analogue of classical lie algebras.
\newblock {\em J. Algebra}, 165:295–--345, 1994.

\bibitem[Lec05]{Lec05}
C{\'e}dric Lecouvey.
\newblock Kostka-foulkes polynomials cyclage graphs and charge statistic for the root system $c_n$.
\newblock {\em J. Alg. Comb.}, 21(2):203--240, 2005.

\bibitem[Lee23]{Lee23}
Seung~Jin Lee.
\newblock Crystal structure on king tableaux and semistandard oscillating tableaux.
\newblock {\em Transformation Groups}, 208(1):438--466, 2023.

\bibitem[LL20]{LL20}
C{\'e}dric Lecouvey and Cristian Lenart.
\newblock Combinatorics of generalize exponents.
\newblock {\em International Mathematics Research Notices}, 2020(16):4942–--4992, 2020.

\bibitem[LOS12]{LOS12}
C{\'e}dric Lecouvey, Masato Okado, and Mark Shimozono.
\newblock Affine crystals, one-dimensional sums and parabolic lusztig q-analogues.
\newblock {\em Mathematische Zeitschrift}, 271:819--865, 2012.

\bibitem[LS78]{LS78}
Alain Lascoux and Marcel-Paul Schützenberger.
\newblock Sur une conjecture de h. o. foulkes.
\newblock {\em C. R. Acad. Sci. Paris Sér. A-B}, 286(7):A323--–A324, 1978.

\bibitem[LS81]{LS81}
Alain Lascoux and Marcel-Paul Schützenberger.
\newblock Le monoïde plaxique.
\newblock {\em Quad. Ricerca Sci.}, 109, 1981.

\bibitem[LS07]{LS07}
C{\'e}dric Lecouvey and Mark Shimozono.
\newblock Lusztig’s $q$-analogue of weight multiplicity and one-dimensional sums for affine root systems.
\newblock {\em Adv. Math.}, 208:438–--466, 2007.

\bibitem[LS24]{LS+}
Seung~Jin Lee and Mark Shimozono.
\newblock Catalan functions for type c.
\newblock {\em In preparation}, 2024+.

\bibitem[Lus83]{Lus83}
George Lusztig.
\newblock Singularities, character formulas, and a $q$-analogue of weight multiplicities.
\newblock {\em Asterisque}, 101--102:208--229, 1983.

\bibitem[Lus90]{Lus90}
George Lusztig.
\newblock Canonical bases arising from quantized enveloping algebras.
\newblock {\em J. Amer. Math. Soc.}, 3:447--498, 1990.

\bibitem[NY97]{NY97}
A.~Nakayashiki and Y.~Yamada.
\newblock Kostka polynomials and energy functions in solvable lattice models.
\newblock {\em Selecta Math.}, 3:547–599, 1997.

\bibitem[Shi05]{Shi05}
Mark Shimozono.
\newblock On the x=m=k conjecture.
\newblock {\em arXiv:math/0501353}, 2005.

\bibitem[SZ06]{SZ06}
Mark Shimozono and Mike Zabrocki.
\newblock Deformed universal characters for classical and affine algebras.
\newblock {\em J. of algebra}, 299:33--61, 2006.

\end{thebibliography}

\end{document}